\documentclass[11pt]{article}
\usepackage{geometry}
\usepackage{graphicx}
\usepackage{amsmath}
\usepackage{amssymb}
\usepackage{amsthm}
\usepackage{float}
\usepackage{epstopdf}
\usepackage{color}
\usepackage{colordvi}
\usepackage{algorithm}
\usepackage{algpseudocode}
\usepackage{subcaption}
\usepackage{siunitx}
\usepackage{comment}
\usepackage{xcolor}

\newcommand{\R}{\mathbb{R}}

\renewcommand{\Vec}[1]{{\boldsymbol{#1}}}
\newcommand{\Mat}[1]{{\boldsymbol{#1}}}

\definecolor{magenta}{rgb}{.5,0,.5}
\definecolor{black}{rgb}{1.0,1.0,1.0}
\definecolor{magenta}{rgb}{.1,0,.3}
\definecolor{gruen}{rgb}{0.2,0.5,.5}
\definecolor{light}{rgb}{ 0.992, 0.961,  0.902}
\definecolor{Tan}{rgb}{ 0.992, 0.9,  0.902}

\title{An Efficient Model Order Reduction Scheme for Dynamic Contact in Linear Elasticity}

\author{D. Manvelyan\footnote{Siemens AG, Corporate Technology, Munich , Germany, diana.manvelyan@siemens.com},    
B. Simeon\footnote{University of Kaiserslautern, Mathematics, Germany, simeon@mathematik.uni-kl.de},     
U. Wever\footnote{Siemens AG, Corporate Technology, Munich , Germany, utz.wever@siemens.com}}

\begin{document}
\maketitle
\begin{abstract}
The paper proposes an approach for the efficient model order reduction of dynamic contact problems in linear elasticity. Instead of
the augmented Lagrangian method that is widely used for mechanical contact problems, we prefer here the Linear Complementarity Programming (LCP) method as basic methodology. It has the advantage of resulting in the much smaller dual problem that is associated with the governing variational principle and that
turns out to be beneficial for the model order reduction.
Since the shape of the contact zone  depends strongly on the acting outer forces, the LCP for the Lagrange multipliers has to be solved in each time step. The model order reduction scheme, on the other hand, is applied to the large linear system for the displacements and computed in advance by means of an Arnoldi process.
In terms of computational effort the reduction scheme is very appealing because the contact constraints are fully satisfied while the reduction  acts only on the displacements. As an extension of our approach, we furthermore take up the idea of the Craig-Bampton method in order to distinguish between interior nodes and the nodes in the contact zone. 
A careful performance analysis closes the paper.
\end{abstract}

\vspace{1cm}
\vspace{1cm}

\section{Introduction}
\label{sec:intro}
The model order reduction for problems in structural mechanics gains rapidly in   importance due to its relevance for simulative operation support \cite{Eigner2012}. Such support is possible if the simulation model runs in parallel to the operation and is synchronized by sensor values at discrete time points. The term "digital twin" has become popular for such a model of the real system. In particular, a digital twin allows to monitor the state of a system at any specified position and at any point in time (and thus, also predict the future states). Obvious benefits of digital twins are, e.g., simplified inspection and service planning, lifetime prediction, advanced fault detection as well as control and optimization during operation. 
The challenge, on the other hand, lies in providing small-scale models that are real time capable and that still preserve the physical key properties of the specific application.

There is a vast literature on model order reduction for different kinds of partial differential equations, see \cite{Lohmann2006,Benner2014,Farhat2014,Chaturantabut2009} for some examples. We concentrate here on reduction methods for dynamic contact problems in linear elasticity.  Despite the linearity of the continuum mechanics model for the displacements, such problems are nonlinear in nature due to the unknown moving contact interface, see, e.g., \cite{Fran75,KiOd88,wriggers2004computational} for background on contact mechanics. 
The reduced model should preserve here both
the system dynamics in terms and the shape of the contact area. Note that for linear elastic models without contact efficient  reduction methods exist, among them modal reduction \cite{Benner2014} and Krylov subspace  methods\cite{Bai2002}, and have been integrated into commercial software \cite{nastran,ansys}. 
The linear mechanical contact problem constitutes a variational problem under unilateral constraints, and its discrete counterpart is mostly solved by the augmented Lagrange method \cite{Gill81}. While this method is very flexible with respect to various descriptions of the contact condition, it may be difficult as a basis for model  order reduction. The main issue with contact problems is that the shape of the contact area is not known a priori, and correspondingly the reduced equations 
need to account for this.

In \cite{Chabal2017}, a dynamic problem with a linear contact condition is discussed and both the displacements and the Lagrange multipliers are reduced  using methods such as the singular value decomposition or the non-negative matrix compression in order to find an online solution for a reduced saddle-point problem by means of the Lagrange multiplier method. In contrast, our approach reduces only the primal displacements. In order to compute the projection matrix only the system matrices and the position vector of the applied force are required. This offline procedure is performed only once. A major benefit is the fact that the trajectories of the full system are not required. Moreover, the computed projection matrix can be re-used for a system with a different load that acts at the same position and in the same direction. After the reduced contact problem is defined, we switch over from the space of reduced displacements to the space of the dual Lagrange multipliers, i.e., to the adjoint problem of the original variational problem.
This leads to a Linear Complementary Programming (LCP) problem \cite{LinearComplementarityProblem,Cottle1992,Li2016}. The complexity of this LCP problem scales with the number of constraints, which is usually small compared to the number of displacements and thus, may be solved even in real time. The resulting update of the displacements is again performed in a reduced space. We think that the degree of reduction is optimal in the sense that the shape of the contact interface is preserved.

For a further separation of the displacement degrees of freedom into interior and contact nodes we additionally apply the method of Craig-Bampton \cite{craigbampton}. It turns out that this separation provides advantages in accuracy in some cases. 
The paper is organized as follows: 
In Section \ref{sec:models} we recall the underlying equations of dynamic contact in linear elasticity. Section \ref{sec:modellcp} describes the solution of the contact problem by Linear Complementary Programming (LCP) both in the static and the dynamic case. Section \ref{sec:reduction} describes the reduction of the problem by Krylov methods and the extension for the Craig-Bampton method. The performance of the approach is demonstrated on 2D geometries in Section \ref{sec:applications}.

\section{Dynamic Contact in Linear Elasticity} 
\label{sec:models}
In this section we summarize the governing equations for dynamic contact
in linear elasticity. We first define a generic model problem in semi-discretized
form that applies to several different finite element discretizations and
then provide some exemplary background.

\subsection{The Generic Model}
We consider large-scale systems of differential-algebraic equations with unilateral constraints in the form  
\begin{eqnarray}\label{eq:KKT_transient1}
&&\Mat{M}\ddot{\Vec{q}}(t) + \Mat{K} \Vec{q}(t)  = \Vec{f}(t)  + \Mat{C}^T\Vec{\lambda}, \\ \label{eq:KKT_transient2}
&&\Mat{C} \Vec{q}(t)  + \Vec{b}\geq 0, \quad \Vec{\lambda} \geq 0, 
	\quad \Vec{\lambda} * (\Mat{C}\Vec{u}(t) +\Vec{b})=0
\end{eqnarray}
where $\Vec{q}(t) \in\R^N$ denotes the vector of nodal displacement variables, 
 $\Mat{M} \in \R^{N\times N}$ the symmetric positive definite mass matrix,
 $\Mat{K} \in \R^{N\times N}$ the symmetric positive semi-definite stiffnes matrix,
 and $\Vec{f}(t) \in \R^N$ the load vector. Furthermore, 
 $\Mat{C}\in \R^{m\times N}$ stands for the constraint matrix and
 $\Vec{b} \in \R^m$ for a given initial clearance vector. The inequality signs
 in (\ref{eq:KKT_transient2}) are to be interpreted component-wise and
 stand for the discretized non-penetration conditions of all variables in the contact interface. Correspondingly, the vector of Lagrange multipliers $\Vec{\lambda} \in \R^m$ 
 enforces the constraint in the dynamic equation (\ref{eq:KKT_transient1}).
 The multiplier $\Vec{\lambda}$, standing for the discretized contact pressure, is always positive and its inner product with the constraints satifies a complementarity condition.
 
In this paper, we will introduce a model order reduction method for the semi-discretized equations (\ref{eq:KKT_transient1}) that projects
the nodal variables $\Vec{q}$ onto a much smaller space while explicitly preserving
the constraints and the Lagrange multipliers. For this purpose, the structure of the equations and the properties of the matrices as stated above are an appropriate starting point. Frictionless, adhesive-free normal contact in combination with small deformation theory and a linear-elastic material will lead to 
(\ref{eq:KKT_transient1}) if 
\begin{itemize}
\item[(i)] the Lagrange multiplier method is used to enforce the non-penetration condition and
\item[(ii)] the finite element method on matching meshes, i.e., node-to-node contact, is applied for the discretization in space.
\end{itemize}
Approaches like the augmented Lagrange method or the Nitsche method 
\cite{wriggers2004computational} lead to modifications 
that will not be considered here. The same holds for additional friction effects and
node-to-segment contact.

For a contact problem with $k$ bodies that fits into the framework described so far,
the mass matrix $\Vec{M}$ will consist of $k$ blocks on the diagonal that  stem from
the discretizations of the individual bodies, and the same block-diagonal structure applies to the stiffness matrix $\Vec{K}$. The constraint matrix $\Vec{C}$ 
will then typically exhibit a sparse structure where each row stands for a node-to-node contact condition.
We do not dive further into the details of various contact models and discretization schemes and refer instead to the standard references \cite{KiOd88,wriggers2004computational}. But for the sake of a self-contained 
presentation, we shortly sketch the setting of the classical obstacle problem
in the dynamic case.

\subsection{Background: Dynamic Obstacle Problem}

Assuming as before
 frictionless, adhesive-free normal contact and a linear-elastic material as well as linear kinematics, we can express the dynamic contact problem in terms 
 of the displacement field $\Vec{u}(\Vec{x},t)\in\R^d$ and the contact pressure 
 $p(\Vec{x},t)\in \R$ where $x \in \Omega \subset \R^d, d=2$ or $d=3$, is the spatial variable 
 and $t \in [t_0,T]$ stands for the temporal variable.
 The boundary of the elastic body is decomposed into
$\partial\Omega = \Gamma_D\cup\Gamma_N\cup\Gamma_C$ where the  
latter represents the contact interface. Using the outward normal vector $\Vec{n}(\Vec{x})$ on
the contact interface, the distance between the body and a given surface is  described by the scalar gap function 
\begin{equation}\label{eq:defgap}
g(\Vec{x}, \Vec{u}) = \big(\Vec{\xi}(\Vec{x}) - (\Vec{x} + \Vec{u}(\Vec{x},t)\big)^T \Vec{n}(\Vec{x}), \quad \Vec{x} \in \Gamma_C.
\end{equation}
Here, the point $\Vec{\xi}(\Vec{x})$ on the obstacle's surface is obtained by projection of $\Vec{x}$ in outward normal direction.
 The contact problem in strong form is then given by the dynamic equations
\begin{equation}\label{elast_cont} 
\rho\Mat{u}_{tt}-\mbox{div} \Mat{\sigma}(\Vec{u})= \Vec{F} \quad \mbox{in } \Omega \times [t_0,T]
\end{equation}
subject to standard boundary conditions
\begin{equation}\label{bcDN}
\Mat{\sigma}(\Vec{u})\cdot\Vec{n}=\Vec{\tau} \quad \mbox{on } \Gamma_N \times [t_0,T], \qquad
\Vec{u}=0  \quad \mbox{on } \Gamma_D \times [t_0,T]
\end{equation}
and furthermore subject to the contact conditions 
\begin{equation}\label{bcContact} 
g \geq 0, \quad
p \geq 0, \quad g \cdot p = 0  \quad \mbox{on } \Gamma_C \times [t_0,T].
\end{equation}
As initial data, 
\begin{equation}
\Mat{u}(\cdot,t_0) = \Mat{u}_0 \quad \mbox{and } \,
\Mat{u}_t(\cdot,t_0) = \Mat{v}_0
\end{equation}
are prescribed with given functions $\Vec{u}_0$ and $\Vec{v}_0$, respectively. 

In (\ref{elast_cont}) and (\ref{bcDN}), the constant $\rho$ denotes the mass density,
$\Vec{F}(\Vec{x},t) \in \R^d$ the volume force,
$\Vec{\tau}(\Vec{x},t) \in \R^d$ the surface traction, 
and 
$\Mat{\sigma}(\Vec{u}) \in R^{d\times d}$ the stress tensor given by
\begin{equation}\label{stress}
  \Mat{\sigma}(\Vec{u}) = \frac{E}{1+\nu} \Mat{e}(\Vec{u}) + \frac{\nu E}{(1+\nu)(1-2 \nu)} \mbox{trace}(\Mat{e}(\Vec{u})) \Mat{I}
\end{equation}
with Young's modulus $E \geq 0$, Poisson's ratio $-1\leq \nu \leq 0.5$ and linearized
strain tensor 
\begin{equation} \label{strain}
  \Mat{e}(\Vec{u})=\frac{1}{2}(\nabla \Vec{u}+\nabla \Vec{u}^T) \in \R^{d\times d}. 
\end{equation}

In order to pass to a weak formulation of the obstacle problem, we introduce
the function spaces
\begin{eqnarray}
\mathcal{V} &:= & \left\{ \Vec{v} \in H^1(\Omega)^d :
                           \Vec{v} = 0 \, \mbox{on } \Gamma_D \right\}, \label{defV} \\
\mathcal{L} &:= & \left\{ \mu \in H^{1/2}(\Gamma_C)^\prime : 
                         \int_{\Gamma_C} \mu w \, ds \geq 0 \,\, \forall
                          w \in H^{1/2}(\Gamma_C), w \geq 0 \right\}
                          \label{defL}
\end{eqnarray}
and the abstract notation
\begin{equation}\label{defANot}
\langle \rho \Vec{u}_{tt}, \Vec{v} \rangle := \int_\Omega \Vec{v}^T \rho \Vec{u}_{tt} dx, \,\,
a(\Vec{u}, \Vec{v}) := \int_\Omega\Mat{\sigma}(\Vec{u}):\Mat{e}(\Vec{v})dx,
\,\,
\langle \ell , \Vec{v} \rangle := \int_\Omega \Vec{v}^T \Mat{F}dx + \int_{\Gamma_N} \Vec{v}^T\Vec{\tau}ds.
\end{equation}
The non-penetration condition in weak form with test function $\mu \in {\cal L}$ is recast as 
\begin{equation}\label{contactweak}
 0 \leq \int_{\Gamma_C} g \mu ds = 
        \int_{\Gamma_C} g (\mu-p) ds  
   = -\int_{\Gamma_C} \Vec{u}^T\Vec{n} (\mu-p) ds
          + \int_{\Gamma_C} (\Vec{\xi}-\Vec{x})^T\Vec{n} (\mu-p) ds.
\end{equation}
The last two integrals give rise to the definitions of the bilinear form
\begin{equation}\label{defBilb}    b_c(\Vec{v},p) :=  - \int_{\Gamma_C}  \Vec{v}^T \Vec{n} \cdot p \,  ds
\end{equation}
on ${\cal V} \times {\cal L}$
and the linear form on $\cal L$ 
\begin{equation}\label{defLFm}
\langle {m}, p \rangle := \int_{\Gamma_C} (\Vec{\xi}-\Vec{x})^T\Vec{n} p \, ds.
\end{equation}
Using these definitions, the weak form of the dynamic contact problem is stated as follows: For each $ t \in [t_0, T] $ find the displacement field $ \Vec{u}(\cdot,t) \in \mathcal{V} $ and the contact pressure $ p(\cdot,t) \in \mathcal{L} $ such that  
\begin{equation}\label{elast_weak}
\begin{aligned} 
\langle \rho \Vec{u}_{tt}, \Vec{v} \rangle + a(\Vec{u}, \Vec{v}) & = \langle \ell,\Mat{v} \rangle  + b_c(\Vec{v}, p) \quad &&\mbox{for all } \Vec{v} \in \mathcal{V},\\
b_c(\Vec{u}, \mu-p) +  \langle {m}, \mu-p \rangle & \geq 0 \quad &&\mbox{for all }  \mu \in \mathcal{L}. 
\end{aligned}
\end{equation}
The system \eqref{elast_weak} is discretized with respect to the spatial varibale by applying the standard Galerkin projection $\Vec{u}(\Vec{x},t) \doteq \sum_{i=1}^N \Vec{\phi}_i(\Vec{x}) q_i(t)
= \Mat{\Phi}(\Vec{x}) \Vec{q}(t)$ with basis functions $\Vec{\phi}_i$ and
nodal variables $\Vec{q} = (q_1, \ldots, q_N)$ to the displacement field.
If the contact condition is simply expressed in terms of the distance between a node
$\Vec{x}_j \in \Gamma_C, i=1,\ldots,m$, and the obstacle, the integral over $\Gamma_C$ 
in the weak dynamic equation
is replaced by a sum
\begin{equation}\label{discretGammaC}
     \int_{\Gamma_C}  \Vec{v}^T \Vec{n} \cdot p \,  ds \doteq
                     \sum_{j=1}^{m}  A_j \Vec{v}(\Vec{x}_j)^T \Vec{n}(\Vec{x}_j) \cdot p(\Vec{x}_j,t), \quad A_j: \mbox{area around node }\Vec{x}_j.
\end{equation}
In the same way, the unilateral constraint is discretized.
Putting finally $\Vec{\lambda}(t) := (p(\Vec{x}_1,t), \ldots, p(\Vec{x}_{m},t))$ as discrete pressure variable,
the semi-discretized equations (\ref{eq:KKT_transient1}) and (\ref{eq:KKT_transient2})
follow in the usual way, with the matrix $\Mat{C} \in \R^{m \times N}$ being an indicator matrix for the nodes in the contact interface and $\Vec{b}\in\R^m$ the offset of the contact. Note that only a few displacements are involved for the contact condition and hence most of the entries of the matrix $\Mat{C}$ are zero.

\section{Solving the Contact Model with LCP} 
\label{sec:modellcp}
A popular approach to solve the mechanical contact problem is the augmented Lagrange method \cite{Powell69,Gill81,wriggers2004computational}, which is very flexible to different descriptions of the contact interface.
However, in the context of model order reduction we advocate a formulation as 
linear complementarity programming (LCP) problem with corresponding solution methods.
We treat first the static case and proceed then to
the dynamic problem.
\subsection{LCP for the Static Contact Problem}
In the stationary case, we can recast the generic semi-discretized formulation  (\ref{eq:KKT_transient1}) and (\ref{eq:KKT_transient2}) as minimization
problem
\begin{equation}\label{eq:var}
\min_{\Vec{q}\in\R^N,\Vec{\lambda}\in\R_+^m} \frac{1}{2}\Vec{q}^T \Mat{K} \Vec{q} - 
\Vec{q}^T \Vec{f} - \Vec{\lambda}^T (\Mat{C}\Vec{q}+\Vec{b}).
\end{equation}
Note that the requirement $\Vec{\lambda}\in\R_+^m$ enforces the positivity of
the contact pressure and, simultaneously, the non-penetration condition.
The corresponding KKT-conditions read
\begin{equation} \label{eq:KKT}
	\begin{array}{l}
	\Mat{K} \Vec{q} - \Vec{f} - \Mat{C}^T \Vec{\lambda}=0, \\
	\Mat{C}\Vec{q}+\Vec{b}\geq 0, \quad \Vec{\lambda} \geq 0, 
	\quad \Vec{\lambda} * (\Mat{C}\Vec{q}+\Vec{b})=0.
	\end{array}
\end{equation}
Next, assuming a positive-definite stiffness matrix, we eliminate the displacements $\Vec{q}$ from the KKT-conditions (\ref{eq:KKT}) via
\begin{equation}\label{eq:DefU}
\Vec{q} = \Mat{K}^{-1} (\Vec{f}  + \Mat{C}^T \Vec{\lambda})
\end{equation}
and insert this expression into the constraint equations. This results in
\begin{equation} \label{eq:KKTL}
	\begin{array}{rcc}
	\Mat{C} \Mat{K}^{-1} \Vec{f} + \Mat{C} \Mat{K}^{-1} \Mat{C}^T \Vec{\lambda}+\Vec{b} & \geq & 0, \\
	\Vec{\lambda} & \geq & 0, \\
	\Vec{\lambda} * (\Mat{C} \Mat{K}^{-1} \Vec{f} + \Mat{C} \Mat{K}^{-1} \Mat{C}^T \Vec{\lambda} + \Vec{b}) & = & 0.
	\end{array}
\end{equation}
Using the abbreviations 
\begin{equation}\label{eq:abbAB}
\Mat{A} := \Mat{C} \Mat{K}^{-1} \Mat{C}^T \in \R^{m \times m}, \quad \Vec{B} := \Mat{C} \Mat{K}^{-1} \Vec{f}+\Vec{b} \in \R^{m},
\end{equation}
the equations (\ref{eq:KKTL}) read
\begin{equation} \label{eq:KKTAB}
	\begin{array}{rcc}
	\Vec{B} + \Mat{A} \Vec{\lambda} & \geq & 0, \\
	\Vec{\lambda} & \geq & 0, \\
	\Vec{\lambda} * (\Vec{B} + \Mat{A} \Vec{\lambda}) & = & 0.
	\end{array}
\end{equation}
 The LCP problem (\ref{eq:KKTAB})  may be solved by standard methods from
 constrained optimization, see below for more details. Note that tasks like the detection of active contact nodes are thus transferred to the LCP solver.

\subsection{LCP for the transient contact problem}

In the transient case, the dynamic equations (\ref{eq:KKT_transient1}) must be discretized in time. As straightforward and unconditionally stable method, we apply the implicit 
Euler scheme, i.e., the first order backward differentiation formula. 
For this purpose, the second order derivative is replaced by the finite difference
\begin{equation}\label{eq:findiff}
\ddot{\Vec{q}}(t+h) \approx \frac{1}{h^2}\big(\Vec{q}(t+h) - 2\Vec{q}(t) + \Vec{q}(t-h)\big)
\end{equation}
with time stepsize $h$. Note that the finite difference $\dot{\Vec{q}}(t+h) \approx
(\Vec{q}(t+h) - \Vec{q}(t))/h$ for the velocity is hidden in (\ref{eq:findiff}), and due to
the second order time derivative, the implicit Euler leads here to a two-step method.
Since (\ref{eq:findiff}) possesses first order of accuracy only, 
one can apply other integration schemes instead, e.g. the generalized-$\alpha $ method, which is very common in structural mechanics, cf.~\cite{simeon}.

In this paper, however, we continue working with the implicit Euler method as our main interest focuses on the  model order reduction scheme. Moreover, dynamic contact involves frequent discontinuities that might impair the convergence of a 
higher order time integration method.
Inserting (\ref{eq:findiff}) into (\ref{eq:KKT_transient1}) leads to
\begin{equation}\label{eq:system_disc}
\Mat{M}\left(\Vec{q}(t+h) -2\Vec{q}(t) + \Vec{q}(t-h)\right) + h^2\Mat{K} \Vec{q}(t+h) = 
h^2\Vec{f}(t+h) + h^2\Mat{C}^T\Vec{\lambda}(t+h)
\end{equation}
Assuming that the previous time steps are known,  (\ref{eq:system_disc}) may be resolved for $\Vec{q}(t+h)$,
\begin{equation}\label{eq:system_disc_resolved}
\Vec{q}(t+h) = (\Mat{M}+h^2\Mat{K})^{-1}(h^2\Vec{f}(t+h) + h^2\Mat{C}^T\Vec{\lambda}(t+h) +2\Mat{M}\Vec{q}(t) - \Mat{M}\Vec{q}(t-h)).
\end{equation}
Now the same procedure as in the static case may be applied. Inserting (\ref{eq:system_disc_resolved}) into the constraints (\ref{eq:KKT_transient2}) again leads to an LCP. With the definitions
\begin{eqnarray}
\Mat{A} &:=& h^2\Mat{C}(\Mat{M}+h^2\Mat{K})^{-1}\Mat{C}^T,\\
\Mat{B} &:=& \Mat{C}(\Mat{M}+h^2\Mat{K})^{-1}(h^2\Vec{f}(t+h) + 
        2\Mat{M}\Vec{u}(t) - \Mat{M}\Vec{u}(t-h)) + \Vec{b}
\end{eqnarray}
we again arrive at the LCP (\ref{eq:KKTAB}), which here has to be solved in each time step.

We remark that the two-step time discretization requires initial values
$\Vec{q}(t_0)$ and $\Vec{q}(t_0+h)$ to start. E.g., if $\Vec{q}_0$ and
$\dot{\Vec{q}}_0$ as initial displacement and velocity are given, one can
compute $\Vec{q}(t_0+h) = \Vec{q}_0 + h \dot{\Vec{q}}_0$ by an explicit Euler step and then continue with the two-step formula (\ref{eq:system_disc}).
The Lagrange multiplier $\Vec{\lambda}$, on the other hand, does not require 
an initial value, and it is computed in each time step as an implicitly given function
of $\Vec{q}$. In the terminology of differential-algebraic equations, this means that 
$\Vec{q}$ stands for the differential variables while the algebraic variables 
$\Vec{\lambda}$ possess no memory. If the unilateral constraints were replaced by 
equality constraints $\Mat{C}\Vec{q} + \Vec{b} = 0$, the index of the resulting
differential-algebraic equation would equal 3, which means that special care 
must be taken for the time integration \cite{BrCP96,HaLR89,simeon}.

Another remark concerns the treatment of impact situations. In the presented
time discretization, impacts are not specifically identified by means of 
computing an impact velocity and a corresponding coefficient of restitution.
Instead, the LCP solution in each time step implicitly enforces a plastic impact 
solution that might lead to a loss in kinetic energy.
\subsection{Remarks on Solving LCP Problems}
The core of our proposed solution algorithm consists of the efficient solution of the LCP problem (\ref{eq:KKTAB}).
One way for solving LCPs is the usage of so-called NCP-functions $\phi^{\text{NCP}}$. An NCP-function $\phi^{\text{NCP}}: \R^2 \to \R$ is characterized by the property
\begin{equation}
\phi^{\text{NCP}}(a,b) = 0 \quad \Leftrightarrow \quad a \geq 0, \ b \geq 0, \ ab = 0.
\end{equation}
Two examples of NCP-functions are:
\begin{eqnarray}
\label{eq_ncp1} \phi_1^{\text{NCP}}(a,b) & = & \min(a,b), \\
\label{eq_ncp2} \phi_2^{\text{NCP}}(a,b) & = & \sqrt{a^2 + b^2} - a - b.
\end{eqnarray}
The first function (\ref{eq_ncp1}) is not suitable for computational purposes because of its non-smoothness. The second NCP-function is also known as the Fischer-Burmeister-function \cite{Fischer92}. Thus for solving an LCP the following two formulations (\ref{eq:solveLCP1}) and (\ref{eq:solveLCP2}) are equivalent:
\begin{equation} \label{eq:solveLCP1}
	\begin{array}{rcc}
	\Vec{B} + \Mat{A} \Vec{\lambda} & \geq & 0, \\
	\Vec{\lambda} & \geq & 0, \\
	\Vec{\lambda} * (\Vec{B} + \Mat{A} \Vec{\lambda}) & = & 0.
	\end{array}
\end{equation}
{\centering $\Longleftrightarrow$}
\begin{equation} \label{eq:solveLCP2}
\phi_2^{\text{NCP}}\left(\Vec{B} + \Mat{A} \Vec{\lambda},\Vec{\lambda}\right) =0.
\end{equation}
Equation (\ref{eq:solveLCP2}) may be solved by Newton's method, but
here we apply a much more efficient method in terms of 
Lemke's algorithm (see \cite{LinearComplementarityProblem}, Sect. 4.4.5).
Specifically, we use the Python implementation of the algorithm \cite{lemke}.

\section{Reduction of the Contact Problem} 
\label{sec:reduction}
In many applications, the number $N$ of degrees of freedom of the discretized system (\ref{eq:KKT_transient1}) is large and its numerical integration is not possible in real time. Model order reduction strategies introduce a reduced state $\Vec{v}\in\R^{n_r}$ with $n_r\ll N$, where $\Vec{v}$ is defined by
\begin{equation}\label{eq:reduction}
\Vec{q} = \Mat{Q}\Vec{v},\quad \Mat{Q}\in\R^{N\times n_r}.
\end{equation}
We can interprete the map from $\Vec{q}$ to $\Vec{v}$ as a second
Galerkin projection that acts on top of the first projection from the 
displacement field $\Vec{u}$ 
to $\Vec{q}$. For a survey on model order reduction see, eg.,  \cite{Baur2014}.

\subsection{Reduction Scheme}
One way to obtain the reduction matrix $\Mat{Q}$ is to use modal reduction. Setting up the eigenvalue problem in the unconstrained case
\begin{equation}
 \omega^2 \Mat{M} \Vec{y}=\Mat{K} \Vec{y} \label{eqn:ev problem}
\end{equation} 
and taking the first $n_r$ eigenvectors, the matrix $\Mat{Q}$ is then defined by
\begin{equation} \label{eq:Beigen}
\Mat{Q}=\left\{\Vec{y}_1,\cdots,\Vec{y}_{n_r} \right\}.
\end{equation} 
A more recent technique are the Krylov subspace methods \cite{Bai2002,Salim2005}. The reduction is defined by
\begin{equation}\label{eq:Bkrylov}
\Mat{Q}=\left\{\Mat{K}_{\omega}^{-1}\Vec{f}, 
\Mat{K}_{\omega}^{-1}\Mat{M} \Mat{K}_{\omega}^{-1}\Vec{f},\cdots,
(\Mat{K}_{\omega}^{-1}\Mat{M})^{n_r-1} \Mat{K}_{\omega}^{-1}\Vec{f} \right\}
\end{equation}
The Krylov vectors can be computed by the Arnoldi algorithm, which delivers an orthonormal basis of the subspace.
Inserting the reduction (\ref{eq:reduction}), where $\Mat{Q}$ is defined by (\ref{eq:Beigen}) or (\ref{eq:Bkrylov}), into the differential equation (\ref{eq:KKT_transient1}) and multiplying by $\Mat{Q}^T$ from the left, we obtain the reduced  equation
\begin{equation}\label{elast_disc_dyn_red}
\hat{\Mat{M}} \ddot{\Vec{v}} + \hat{\Mat{K}}\Vec{v} =
\Mat{Q}^T\left(\Vec{f}(t)+\Mat{C}^T \Vec{\lambda}(t)\right)
\end{equation}
where $\hat{\Mat{M}}=\Mat{Q}^T \Mat{M}\Mat{Q}$, $\hat{\Mat{K}}=\Mat{Q}^T \Mat{K}\Mat{Q}$ $\in \R^{n_r\times n_r}$ are the reduced  matrices.  

Using (\ref{elast_disc_dyn_red}),  moreover the reduced time discretization 
follows easily as 
\begin{equation}\label{eq:system_time_red}
\Vec{v}(t+h) = \left(\hat{\Mat{M}}+h^2\hat{\Mat{K}}\right)^{-1}
\left(h^2 \Mat{Q}^T\Vec{f}(t+h) + h^2\Mat{Q}^T\Mat{C}^T\Vec{\lambda}(t+h) +
 2\hat{\Mat{M}}\Vec{v}(t) - \hat{\Mat{M}}\Vec{v}(t-h))\right)
\end{equation}
We stress that (\ref{eq:system_time_red}) contains only operations with complexity $O(n_r)$ or $O(m)$.
Though $\Mat{Q}\in\R^{N\times n_r}$, the acting force $\Vec{f}(t+h)$  has only a few entries and the Lagrange multiplier $ \Vec{\lambda}$ is still of dimension $m$.

\subsection{Reduction Scheme and Craig-Bampton Method} \label{craig_bampton}
Next, we introduce a modification of our approach 
which is inspired by the fact that often the dominant influence on the dynamics of the system stems from the nodes located on the contact interface. In other words, we distinguish between the nodal variables on the contact interface, so-called master nodes and the ones outside of the contact interface, the so-called slave nodes. 
We adopt this technique from structural dynamics, see Craig and Bampton \cite{craigbampton}.

Let $ \Vec{q} = (\Vec{q}_M, \Vec{q}_S)^T \in \mathbb{R}^N $ the permuted solution vector such that we distinguish between master and slave nodes. Then the corresponding permuted matrices and the permuted force vector are given by
\begin{align}
\Mat{M} = 
\begin{pmatrix}
\Mat{M}_{MM} & \Mat{M}_{MS}\\
\Mat{M}_{SM} & \Mat{M}_{SS}
\end{pmatrix},
\quad
\Mat{K} = 
\begin{pmatrix}
\Mat{K}_{MM} & \Mat{K}_{MS}\\
\Mat{K}_{SM} & \Mat{K}_{SS}
\end{pmatrix},
\quad
\Mat{F} =
\begin{pmatrix}
\Mat{F}_M\\
\Mat{F}_S
\end{pmatrix}.
\end{align} 
Due to the fact that the constraint matrix affects only the dynamics of the slave nodes it holds
\begin{align}\label{no_const}
\Mat{C} = 
\begin{pmatrix}
\Mat{C}_M\\
\Mat{C}_S
\end{pmatrix},
\quad
\Mat{C}_S = 0.
\end{align}
The constraint matrix refers only to the master nodes.
A crucial assumption is now that the influence acceleration term for the slave nodes 
in the dynamic equation can be neglected, and the same for the force  term.
 This yields a coupling equation between the master and slave node via
\begin{equation} \label{eq:coupling}
\Mat{K}_{SM} \Vec{q}_{M} + \Mat{K}_{SS} \Vec{q}_{S} = 0
\end{equation}
Since we want to preserve the contact nodes and only reduce the slave nodes, we keep the structure of the master nodes fixed, i.e, $ \Vec{u}_M = 0. $
Then due to \eqref{no_const} the dynamical system is given by
\begin{equation} \label{eq:slave}
\Mat{M}_{SS}\ddot{\Vec{q}}_{S} + \Mat{K}_{SS}\Vec{q}_S = \Mat{F}_S
\end{equation}
 Instead for the full system, we compute the transformation matrix for the slave system \eqref{eq:slave} using the Arnoldi method, which yields
\begin{equation}\label{eq:Bkrylov_S}
\Mat{Q}_S:=\left\{\Mat{K}_{SS}^{-1}\Vec{F}_S, 
\Mat{K}_{SS}^{-1}\Mat{M}_{SS} \Mat{K}_{SS}^{-1}\Vec{F}_S,\cdots,
(\Mat{K}_{SS}^{-1}\Mat{M}_{SS})^{n-1} \Mat{K}_{SS}^{-1}\Vec{F}_S \right\}.
\end{equation}
Using \eqref{eq:coupling}, we obtain the complete transformation matrix 
that includes the relation between the master and slave nodes,
\begin{align}\label{eq:defQCB}
\Mat{Q}_{CB} =
\begin{pmatrix}
\Mat{I}_M &  0 \\
-\Mat{K}_{SS}^{-1}\Mat{K}_{SM} & \Mat{Q}_S
\end{pmatrix}
\end{align}
The transformation matrix $ \Mat{Q}_{CB} $ with the partitioning of master and slave nodes in the fashion of the Craig-Bampton method provides an alternative 
for reducing the transient contact problem \eqref{eq:system_time_red}.
Regardless whether we employ $\Mat{Q}$ from the reduction (\ref{eq:Bkrylov}) or 
$ \Mat{Q}_{CB} $, the computation of the transformation can be done in an a priori offline phase and does not require solution trajectories, i.e., snapshots, of the full system.

\section{Applications}
\label{sec:applications}
In this section we apply our approach to 2D problems with self-contact. The spatial  discretization uses standard bilinear shape functions on quadrilateral elements.
The finite element method  and the reduction scheme are implemented by means of a custom Python script while Lemke's algorithm is taken from an open source library \cite{lemke}. 
\subsection{Simulations of the Full and Reduced Model}
We consider the unit square $\Omega = [0,1] \times [0,1]$ under the assumption of plane stress. The parameters in dimensionless form are $\rho =1,  E = 1000, \nu = 0.3 $. The left edge is fixed by means of zero Dirichlet boundary conditions.
We discretize the example by 1600 quadrilateral elements, which lead to $N = 3386$ DOF. An interior contact interface is predefined by a fixed number of $ m $  discretization points where each such point is represented by double nodes. This kind of data structure allows us to distinguish the nodes placed on the contact interface from the remaining nodes. Note that the node-to-node non-penetration condition acts only in $x$-direction while the DOF in $y$-direction may move freely.

\begin{figure}[h!]
	\centering
	\includegraphics[angle=0,scale=0.5]{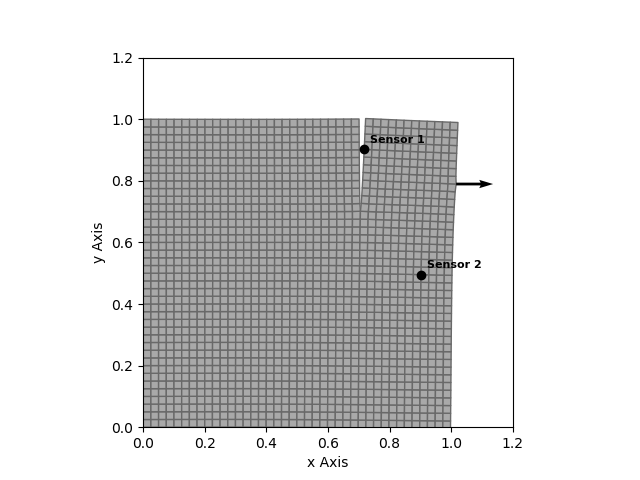}
	\caption{  \label{figgeom1}The depicted deformed state with an open contact tear occurs at the time point when the acting horizontal nodal load vector is pointing outward from the domain.}
\end{figure}

An oscillating force acts on the right side of the domain which affects the contact interface by opening and closing the tear over time.
The problem setup and the mesh are illustrated in Figure \ref{figgeom1}. Two sensors are placed for the data extraction in order to compare the displacements of the full and the reduced model.
Before starting the reduction method we eliminate the fixed degrees of freedom due to the Dirichlet boundary condition. Overall, there are $ 4m $ displacement variables associated with the contact nodes due to the relation $ n_c = 4m$

\subsubsection{\textit{First Example}}

\begin{figure}[h!]
	\centering
	\includegraphics[angle=0,scale=0.4]{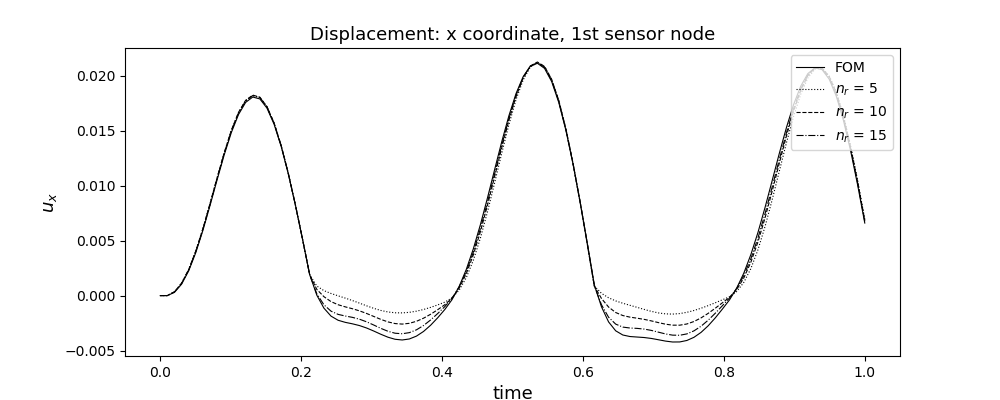}
	\bigbreak
	\includegraphics[angle=0,scale=0.4]{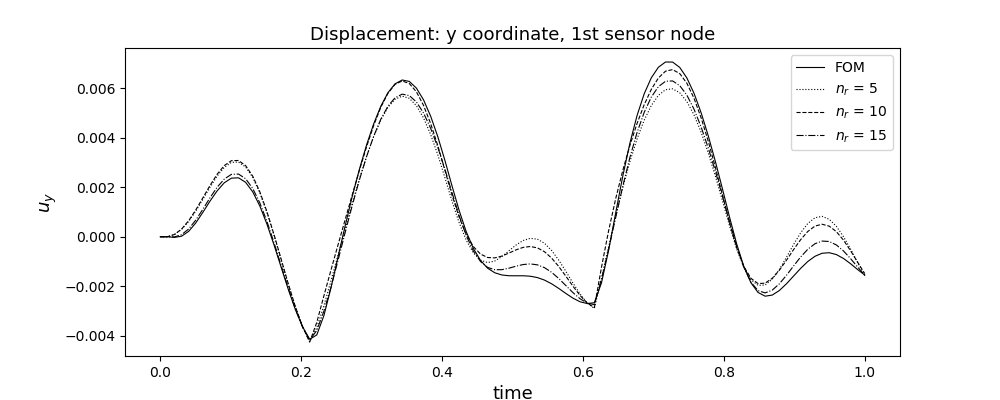}
	\caption{The $x$- and $y$-displacement  of the full order model (FOM) and the reduced order model for $ n_r = 5, 10, 15 $ basis vectors. The displacement trajectory is evaluated in the first sensor node located on the contact interface.  }\label{fig:Krylov1}
\end{figure}

In our first example the contact tear is located on the upper-right side of the domain. The contact is discretized by $ m = 12 $ points and there are $n_c = 48$ contact variables respectivly.
We consider an oscillating nodal force in horizontal direction (cf. Figure \ref{figgeom1}) with the load magitude given by
\begin{equation}\label{load1}
\Vec{f}(t) = \left(\begin{array}{c}
              1.5\sin(0.1\pi t) \\ 0 \end{array}\right), \quad t \in [0,1].
\end{equation}
We test our approach for three different numbers $ n_r = 5, 10, 15$ of the Krylov basis vectors. The resulting trajectories for the 
displacement in the first sensor node in the contact interface, Figure \ref{fig:Krylov1}, show very good agreement of the reduced order
model with the full order model.

\begin{figure}[h!]
	\centering
	\includegraphics[angle=0,scale=0.4]{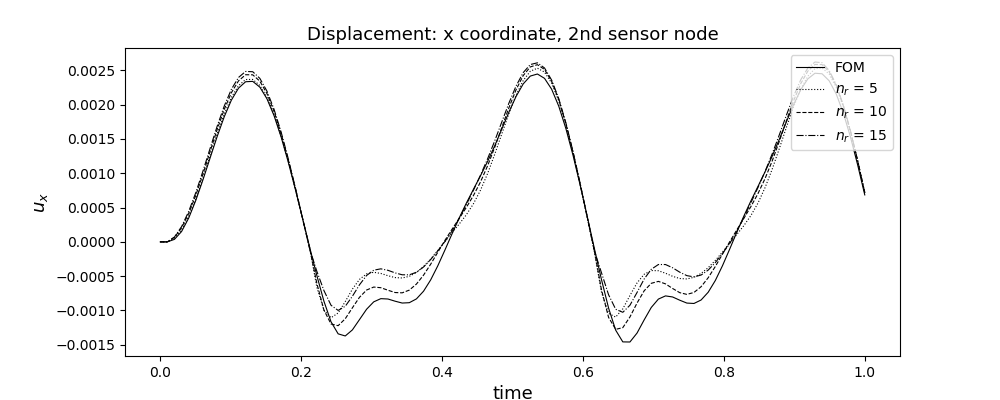}
	\bigbreak
	\includegraphics[angle=0,scale=0.4]{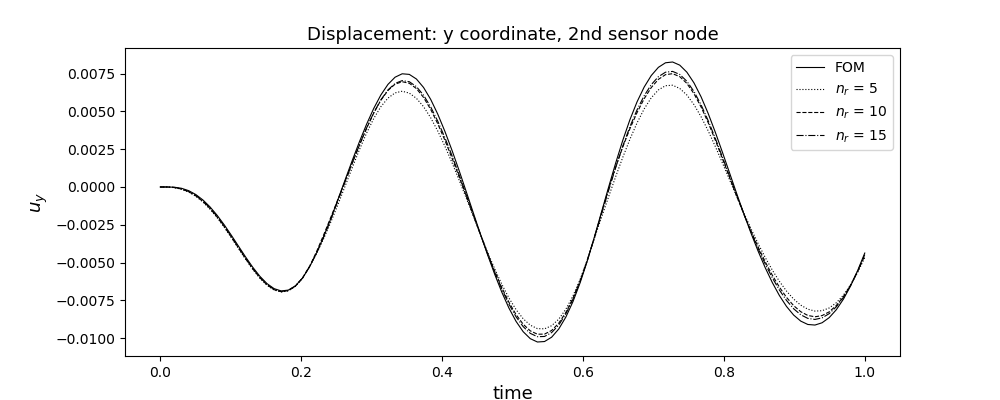}
	\caption{The $x$- and $y$-displacements of the full order model (FOM) and the reduced order model for $ n_r = 5, 10, 15$ basis vectors. The displacement trajectories are evaluated in the second sensor node located on the contact interface.}\label{fig:Krylov2}
	
\end{figure}

As expected, the approximation quality increases with the dimension of the Krylov subspace. Nevertheless, the displacement period seems to be recovered even for small numbers of $ n_r $. The results for the second sensor node are also very satisfactory, as presented in
Fig. \ref{fig:Krylov2}.

\begin{figure}[h!]
	\centering
	\includegraphics[angle=0,scale=0.4]{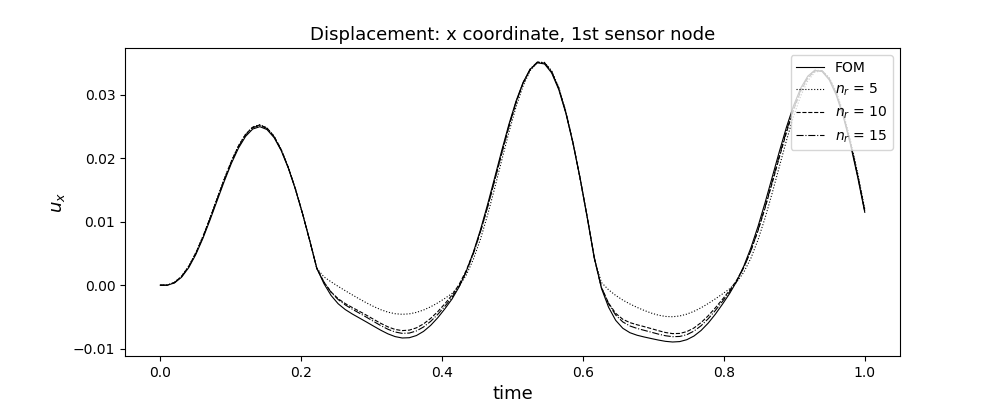}
	\bigbreak
	\includegraphics[angle=0,scale=0.4]{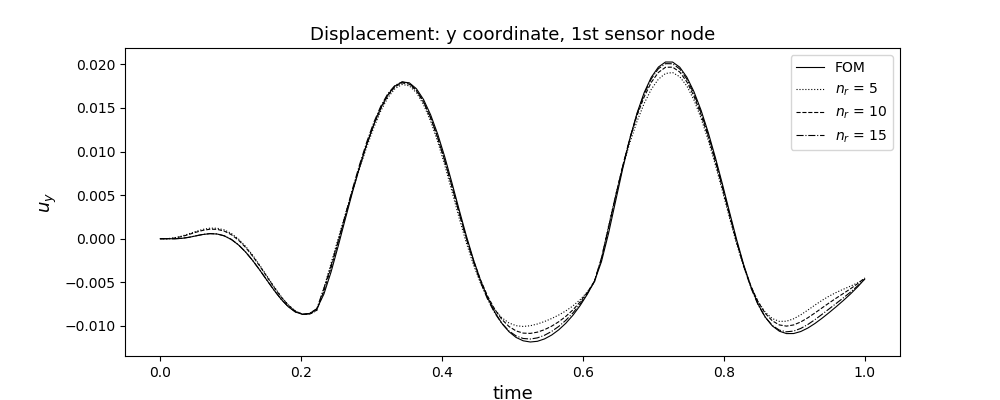}
\caption{$x$- and $y$-displacements for the FOM (full order model) and for the ROM (reduced order model) with three different numbers of basis vectors, i.e., $ n_r = 5, 10, 15.$
First sensor node and load (\ref{load2}).}\label{fig:Krylov3}
\end{figure}

\begin{figure}[h!]
	\centering
	\includegraphics[angle=0,scale=0.4]{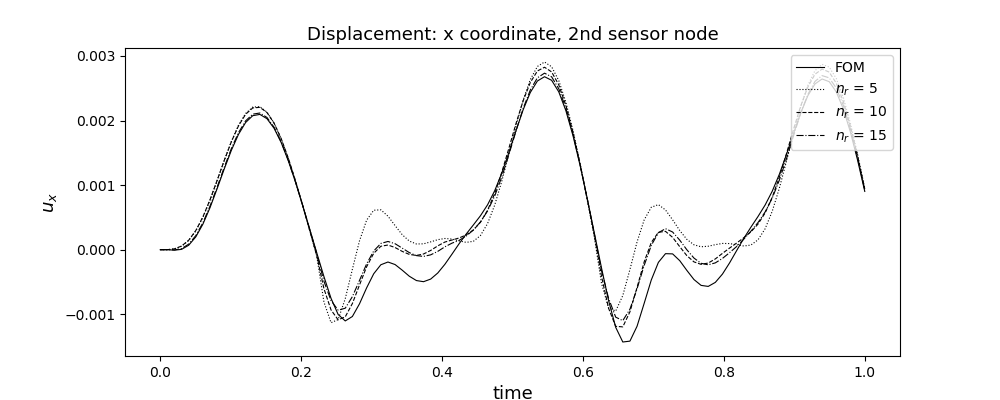}
	\bigbreak
	\includegraphics[angle=0,scale=0.4]{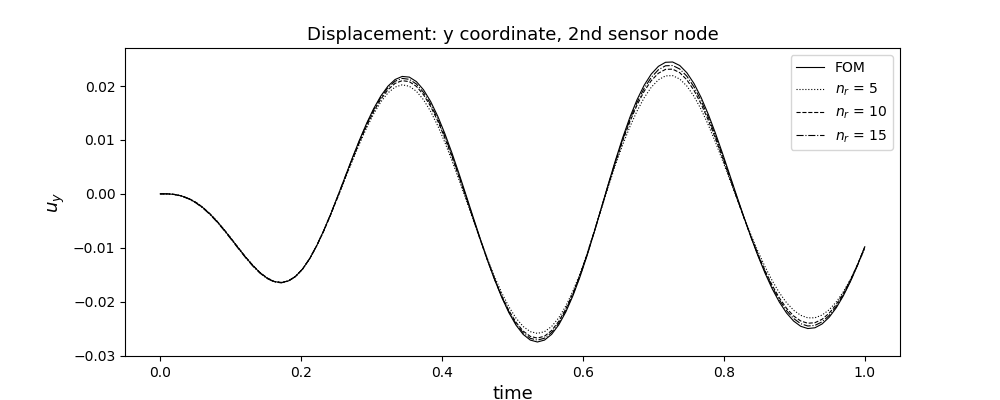}
	\caption{$x$- and $y$-displacements for the FOM (full order model) and for the ROM (reduced order model) with three different numbers of basis vectors, i.e., $ n_r = 5, 10, 15.$
	Second sensor node and load  (\ref{load2})}\label{fig:Krylov4}
\end{figure}

In our tests we also applied an outer force with the same magnitude as \eqref{load1} but different angle given by
\begin{equation}\label{load2}
\Vec{f}(t) = \left(\begin{array}{c}
           1.5\sin(0.1\pi t)\sin(0.25 \pi)\\
            1.5\sin(0.1\pi t)\sin(0.25 \pi) \end{array}\right).
\end{equation}
In this case the $y$-displacement shows a stronger response.
The comparison of the displacements of the FOM and ROM is depicted in Figure \ref{fig:Krylov3}
and \ref{fig:Krylov4} corresponding to the first sensor and to the second sensor, respectively.


\newpage
\subsubsection{\textit{Second Example}}

\begin{figure}[h!]
	\centering
	\includegraphics[angle=0,scale=0.5]{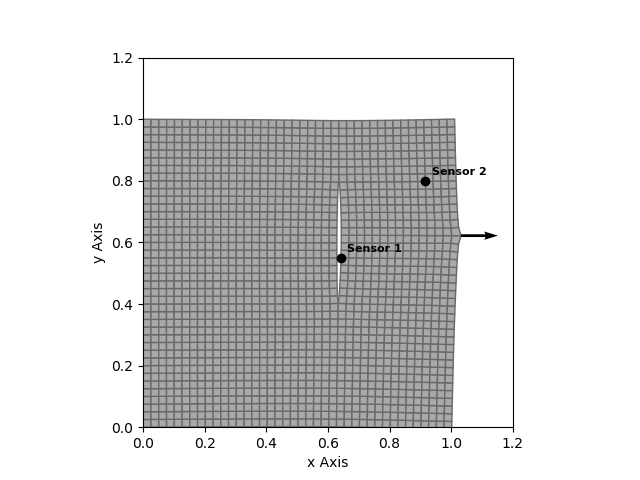}
	\caption{  \label{figgeom3} Simulation scenario with interior tear. The depicted state is captured at the time instant when the horizontal nodal load is directed outside the domain, opening the interior contact tear.}
\end{figure}

\begin{figure}[h!]
	\centering
	\includegraphics[angle=0,scale=0.4]{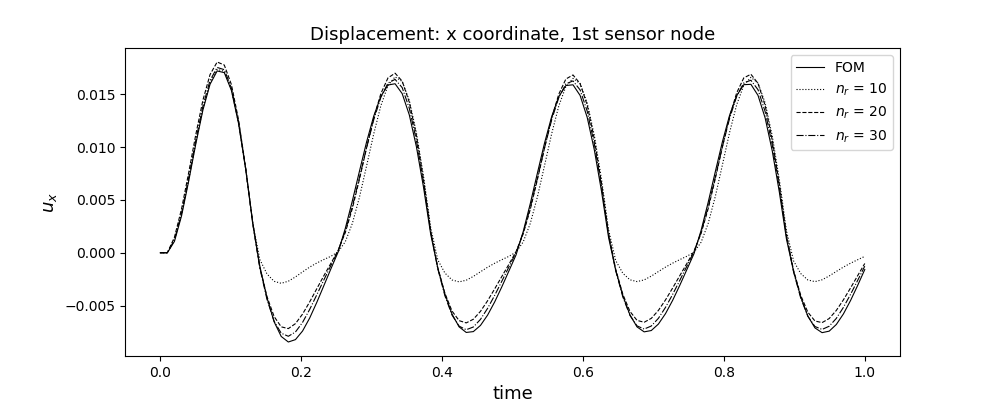}
	\bigbreak
	\includegraphics[angle=0,scale=0.4]{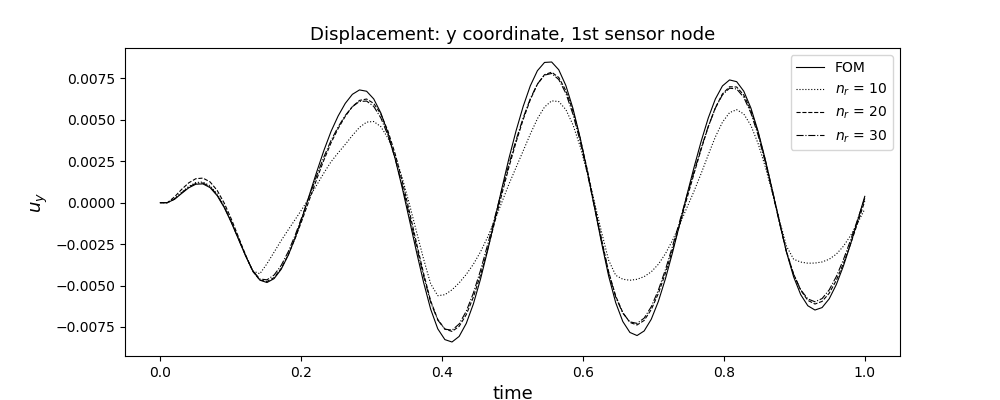}
	\caption{The $x$- and $y$-displacements of the FOM and ROM for $ n_r = 10, 20, 30$ basis vectors respectvely. The displacement trajectories are evaluated on the first sensor node located on the inner contact interface and correspond to the acting force \eqref{load1} and Figure \ref{figgeom3}. }\label{fig:Krylov5}
\end{figure}

Next, we consider another example where the contact interface is placed inside the domain, see Fig. \ref{figgeom3}. The corresponding displacements can be found in Fig. \ref{fig:Krylov5}.  

In contrast to the examples before, the double nodes appear in the interior of the domain. The outer boundary and the surrounding inner area are represented by single nodes.
The stiffness of the surface and connectedness of the nodes produces an inner tear.
Again, the reduced order model performs well. At least $n_r = 20$ basis vectors are 
required to obtain satisfactory results.
\subsection{Splitting up the Projection Matrix: Craig-Bampton}
\begin{figure}[h!]
	\centering
	\includegraphics[angle=0,scale=0.4]{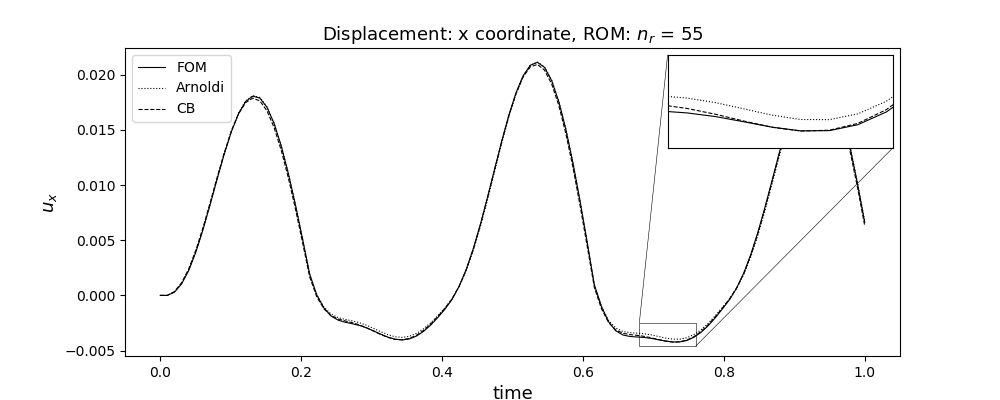}
	\bigbreak
	\includegraphics[angle=0,scale=0.4]{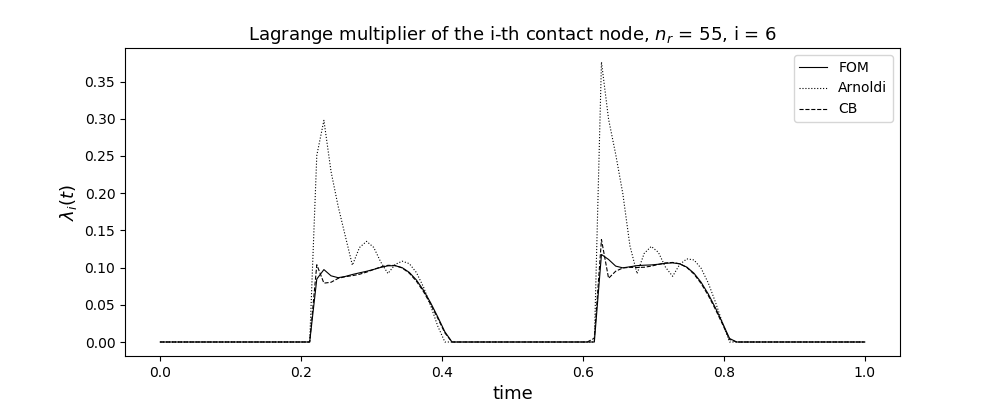}
	\caption{Upper figure: $x$-displacements of the FOM and the ROM, once without (Arnoldi) and once with splitting (CB). In both cases, $ n_r = 55. $ In case of CB $ n_r = n_k + n_c $ with $ n_k = 7 $ and $ n_c = 48.$ Trajectory of first sensor node, which is placed on the contact interface. Lower figure:  Lagrange multiplier in contact node $6$.}\label{fig:Ar_vs_CB1}
\end{figure}

We continue with a comparison of the reduction results with and without Craig-Bampton splitting, which was discussed in  Subsection \ref{craig_bampton}. To guarantee a fair comparison of the two methods, we choose the same number of basis functions in both cases. In case of using no splitting, the basis consists of only Krylov vectors, whereas in the splitting case, the basis consists of a few Krylov vectors ($ n_k $) and additionally the contact nodes ($ n_c $ ), i.e., $ n_r = n_k + n_c. $ The Figure \ref{fig:Ar_vs_CB1} shows the comparison of both reduction approaches with the full model,
referring to the scenario of Fig. \ref{figgeom1}.
We can observe an increase of accuracy in case of using the splitting of Craig-Bampton. 
This holds also for the second sensor node and the displacements in $y$-direction.
 
Furthermore, the Lagrange multipliers $ \Vec{\lambda}(t) $  exhibit an 
interesting behavior. We observe that ROM without splitting produces contact pressures
that are quite different from those of the FOM, while the Craig-Bampton splitting yields values
that agree well with the FOM. Both approaches guarantee positivity of the Lagrange multiplicator whenever the contact is activated. 
An explanation of  this deviation with respect to the Lagrange multipliers in case of no splitting
is as follows. The projection matrix $ \Mat{Q} $ in \eqref{eq:reduction} alters the nodewise distribution of the pressure on the contact variables considerably and involves spatial basis functions with large support.
Splitting of variables, on the other hand, by means of $\Mat{Q}_{CB}$ in
(\ref{eq:defQCB})
preserves the contact variables and thus the Lagrange multipliers. In Figure \ref{fig:Ar_vs_CB1} 
the $x$-displacement in the sensor node 1 for both reduction approaches and the FOM is displayed,
along with the corresponding Lagrange multiplier.
The switching between contact, where the multiplier is positive, and positive gap with 
vanishing multiplier can be nicely observed.

\section{Conclusions}
\label{sec:conclusions}
In this paper we have presented a reduction method for dynamic contact problems
in linear elasticity. 
The generic model elaborated in Section 2 does not cover all available methods but nevertheless is general enough to provide sufficient insight into the problem class.
We plan to extend the problem class in the future to situations like node-to-segment contact and non-matching meshes where the unilateral constraints might become 
nonlinear.

The proposed model order reduction uses a Krylov subspace approach for the displacement variables but preserves the Lagrange multipliers. This is intimately connected with the LCP problem that we advocate for solving the contact problem in each time step. Obviously, the savings in DOF are substantial as long as 
the number of contact constraints is small compared to the number of displacement variables. If we consider solids and contact on outer or interior boundaries, this holds true, but if we extend the method to shell problems where contact might occur all over the computational domain, the savings might be less substantial.

It also turned out that the numerical approximation of the Lagrange multipliers, i.e., the contact pressure in the contact interface, is in general not accurate for
the simple Krylov reduction scheme. By splitting of the nodes in the fashion of the Craig-Bampton method, a substantial improvement can be achieved. Moreover,
this modification leads also to better approximations of the displacement variables.
In any case, the reduction matrix is computed only once in an offline procedure.
There is no need to compute trajectories or snapshots of the full model to provide
data for the reduction scheme.

Our future effort will focus on the development of methods for treating more complex geometries and on the analysis of the proposed method class in terms of 
apriori and a posteriori error estimates.

\vspace{1.0cm}
{\bf Acknowledgement}: 
The authors would like to thank Meinhard Paffrath for his guidance on the technical matters concerning the implementation. Furthermore, the authors want to thank Christoph Heinrich and his analysis group for valuable discussions.

\vspace{1.0cm};
\bibliographystyle{plain}
\bibliography{references}

\end{document}